\documentclass[reqno]{amsart}
\usepackage[left=1in, right=1in, top=1in]{geometry}

\usepackage{etoolbox}

\makeatletter
\let\old@tocline\@tocline
\let\section@tocline\@tocline
\newcommand{\section@dotsep}{4.5}
\newcommand{\subsection@dotsep}{4.5}
\patchcmd{\@tocline}
  {\hfil}
  {\nobreak
     \leaders\hbox{$\m@th
        \mkern \section@dotsep mu\hbox{.}\mkern \section@dotsep mu$}\hfill
     \nobreak}{}{}
\let\section@tocline\@tocline
\let\@tocline\old@tocline

\patchcmd{\@tocline}
  {\hfil}
  {\nobreak
     \leaders\hbox{$\m@th
        \mkern \subsection@dotsep mu\hbox{.}\mkern \subsection@dotsep mu$}\hfill
     \nobreak}{}{}
\let\subsection@tocline\@tocline
\let\@tocline\old@tocline

\let\old@l@section\l@section
\let\old@l@subsection\l@subsection

\def\@tocwriteb#1#2#3{%
  \begingroup
    \@xp\def\csname #2@tocline\endcsname##1##2##3##4##5##6{%
      \ifnum##1>\c@tocdepth
      \else \sbox\z@{##5\let\indentlabel\@tochangmeasure##6}\fi}%
    \csname l@#2\endcsname{#1{\csname#2name\endcsname}{\@secnumber}{}}%
  \endgroup
  \addcontentsline{toc}{#2}%
    {\protect#1{\csname#2name\endcsname}{\@secnumber}{#3}}}%

\newlength{\@tocsectionindent}
\newlength{\@tocsubsectionindent}
\newlength{\@tocsubsubsectionindent}
\newlength{\@tocsectionnumwidth}
\newlength{\@tocsubsectionnumwidth}
\newlength{\@tocsubsubsectionnumwidth}
\newcommand{\settocsectionnumwidth}[1]{\setlength{\@tocsectionnumwidth}{#1}}
\newcommand{\settocsubsectionnumwidth}[1]{\setlength{\@tocsubsectionnumwidth}{#1}}
\newcommand{\settocsubsubsectionnumwidth}[1]{\setlength{\@tocsubsubsectionnumwidth}{#1}}
\newcommand{\settocsectionindent}[1]{\setlength{\@tocsectionindent}{#1}}
\newcommand{\settocsubsectionindent}[1]{\setlength{\@tocsubsectionindent}{#1}}
\newcommand{\settocsubsubsectionindent}[1]{\setlength{\@tocsubsubsectionindent}{#1}}

\renewcommand{\l@section}{\section@tocline{1}{\@tocsectionvskip}{\@tocsectionindent}{}{\@tocsectionformat}}%
\renewcommand{\l@subsection}{\subsection@tocline{1}{\@tocsubsectionvskip}{\@tocsubsectionindent}{}{\@tocsubsectionformat}}%
\renewcommand{\l@subsubsection}{\subsubsection@tocline{1}{\@tocsubsubsectionvskip}{\@tocsubsubsectionindent}{}{\@tocsubsubsectionformat}}%
\newcommand{\@tocsectionformat}{}
\newcommand{\@tocsubsectionformat}{}
\newcommand{\@tocsubsubsectionformat}{}
\expandafter\def\csname toc@1format\endcsname{\@tocsectionformat}
\expandafter\def\csname toc@2format\endcsname{\@tocsubsectionformat}
\expandafter\def\csname toc@3format\endcsname{\@tocsubsubsectionformat}
\newcommand{\settocsectionformat}[1]{\renewcommand{\@tocsectionformat}{#1}}
\newcommand{\settocsubsectionformat}[1]{\renewcommand{\@tocsubsectionformat}{#1}}
\newcommand{\settocsubsubsectionformat}[1]{\renewcommand{\@tocsubsubsectionformat}{#1}}
\newlength{\@tocsectionvskip}
\newcommand{\settocsectionvskip}[1]{\setlength{\@tocsectionvskip}{#1}}
\newlength{\@tocsubsectionvskip}
\newcommand{\settocsubsectionvskip}[1]{\setlength{\@tocsubsectionvskip}{#1}}
\newlength{\@tocsubsubsectionvskip}
\newcommand{\settocsubsubsectionvskip}[1]{\setlength{\@tocsubsubsectionvskip}{#1}}

\patchcmd{\tocsection}{\indentlabel}{\makebox[\@tocsectionnumwidth][l]}{}{}
\patchcmd{\tocsubsection}{\indentlabel}{\makebox[\@tocsubsectionnumwidth][l]}{}{}
\patchcmd{\tocsubsubsection}{\indentlabel}{\makebox[\@tocsubsubsectionnumwidth][l]}{}{}

\newcommand{\@sectypepnumformat}{}
\renewcommand{\contentsline}[1]{%
  \expandafter\let\expandafter\@sectypepnumformat\csname @toc#1pnumformat\endcsname%
  \csname l@#1\endcsname}
\newcommand{\@tocsectionpnumformat}{}
\newcommand{\@tocsubsectionpnumformat}{}
\newcommand{\@tocsubsubsectionpnumformat}{}
\newcommand{\setsectionpnumformat}[1]{\renewcommand{\@tocsectionpnumformat}{#1}}
\newcommand{\setsubsectionpnumformat}[1]{\renewcommand{\@tocsubsectionpnumformat}{#1}}
\newcommand{\setsubsubsectionpnumformat}[1]{\renewcommand{\@tocsubsubsectionpnumformat}{#1}}
\renewcommand{\@tocpagenum}[1]{%
  \hfill {\mdseries\@sectypepnumformat #1}}

\let\oldappendix\appendix
\renewcommand{\appendix}{%
  \leavevmode\oldappendix%
  \addtocontents{toc}{%
    \protect\settowidth{\protect\@tocsectionnumwidth}{\protect\@tocsectionformat\sectionname\space}%
    \protect\addtolength{\protect\@tocsectionnumwidth}{2em}}%
}
\makeatother

\makeatletter
\settocsectionnumwidth{2em}
\settocsubsectionnumwidth{2.5em}
\settocsubsubsectionnumwidth{3em}
\settocsectionindent{1pc}%
\settocsubsectionindent{\dimexpr\@tocsectionindent+\@tocsectionnumwidth}%
\settocsubsubsectionindent{\dimexpr\@tocsubsectionindent+\@tocsubsectionnumwidth}%
\makeatother

\settocsectionvskip{0pt}
\settocsubsectionvskip{0pt}
\settocsubsubsectionvskip{0pt}

\usepackage{mathrsfs}
\usepackage{mathrsfs}
\usepackage{amsfonts}
\usepackage{comment}
\usepackage{cases}
\usepackage{latexsym}
\usepackage{amsmath}
\usepackage[all]{xy}
\usepackage{stmaryrd}
\usepackage{amsfonts}
\usepackage{amsmath,amssymb,amscd,bbm,amsthm,mathrsfs,dsfont}

\usepackage{tikz}

\usepackage{pgflibraryarrows}
\usepackage{pgflibrarysnakes}

\usepackage[numbers,sort&compress]{natbib}
\usepackage[pagebackref]{hyperref}
\usepackage{hypernat}

\usepackage{color, hyperref}
\definecolor{blue}{rgb}{0.9,0.0,0.9}
\hypersetup{colorlinks, breaklinks,
            linkcolor=red,urlcolor=blue,
            anchorcolor=blue,citecolor=blue}

\usepackage{fancyhdr}
\usepackage{amsxtra,ifthen}
\usepackage{verbatim}

\usepackage{fancyhdr}
\usepackage{amsxtra,ifthen}
\usepackage{verbatim}

\numberwithin{equation}{section}

\theoremstyle{plain}
\newtheorem{theorem}{Theorem}[section]
\newtheorem{lemma}{Lemma}[section]

\newtheorem{corollary}{Corollary}[section]

\theoremstyle{definition}

\allowdisplaybreaks[4]

\begin{document}
\title{Nonlinear skew Lie-type centralizers on $\ast$-algebras}%

\author{Mohammad Afajal Ansari,  Mohammad Ashraf and Md Shamim Akhter}

\address{M. A. Ansari, Department of Mathematics, Shibli National College, Azamgarh 276001, India}
\email{afzalgh1786@gmail.com}

\address{M. Ashraf, Department of Mathematics, Aligarh Muslim University, Aligarh 202002, India}
\email{mashraf80@hotmail.com}

\address{M. S. Akhter, Department of Mathematics, Aligarh Muslim University, Aligarh 202002, India}
\email{akhter2805@gmail.com}

\pagestyle{plain}


\begin{abstract}
Let \(\mathcal{A}\) be a unital \(\ast\)-algebra over the complex field \(\mathbb{C}\), with center \(\mathcal{Z}(\mathcal{A})\). 
For any \(X_1,X_2,\ldots,X_n\in\mathcal{A}\), define recursively
$
p_1(X_1)=X_1,
$
$p_2(X_1, X_2)=[X_1,X_2]_{\ast}=X_1X_2-X_2X_1^{\ast}$
and
$
p_n(X_1,X_2,\ldots,X_n)
=[p_{n-1}(X_1,X_2,\ldots,X_{n-1}),X_n]_{\ast}
$
for all integers \(n\geq 2\).
In this paper, we prove that if a map \(\Phi:\mathcal{A}\to\mathcal{A}\) satisfies
$
\Phi(p_n(X_1,X_2,\ldots,X_n))
=p_n(\Phi(X_1),X_2,\ldots,X_n)
$
for all \(X_1,X_2,\ldots,X_n\in\mathcal{A}\), then \(\Phi\) necessarily has the form
$
\Phi(X)=ZX,
$
where \(Z=Z^{\ast}\in\mathcal{Z}(\mathcal{A})\).
As an application, we establish that, under suitable assumptions, every nonlinear generalized skew Lie \(n\)-derivation \(G:\mathcal{A}\to\mathcal{A}\) admits a standard form
$
G(X)=ZX+D(X),
$
where \(Z\in\mathcal{Z}(\mathcal{A})\) and \(D:\mathcal{A}\to\mathcal{A}\) is an additive \(\ast\)-derivation. Finally, we apply this result to several important classes of unital \(\ast\)-algebras, including unital prime \(\ast\)-algebras, standard operator algebras, and von Neumann algebras.
\end{abstract}

\maketitle

\section{Introduction}
Let \(\mathcal{A}\) be an associative \(\ast\)-algebra over the complex field \(\mathbb{C}\). For any \(X,Y\in\mathcal{A}\), the skew Lie product of \(X\) and \(Y\) is defined by
$
[X,Y]_{\ast}=XY-YX^{\ast}.
$
Several related products induced by the involution \(\ast\), such as the skew Jordan product
$
X\ast Y=XY+YX^{\ast},
$
the bi-skew Lie product
$
[X,Y]_{\bullet}=XY^{\ast}-YX^{\ast},
$
and the bi-skew Jordan product
$
X\bullet Y=XY^{\ast}+YX^{\ast},
$
have been extensively investigated in various classes of \(\ast\)-algebras (see, for instance, \cite{Ashraf2022,Ashraf2023,Ashraf2024,Bresar2000,LiZhao2021,Molnar97,MS97,PS90,PS91,PS93,PS94,YZ12,Zhang16}). These products play a significant role in several areas of mathematical research, and their structural properties have consequently attracted considerable attention. In particular, the skew Lie product has been studied extensively because, as demonstrated by the fundamental theorem of {\v S}emrl \cite{PS91}, mappings of the form \(T\mapsto TA-AT^{\ast}\) arise naturally in the problem of representing quadratic functionals by sesquilinear functionals (see \cite{PS90,PS93}).

Recall that an additive map $D:\mathcal{A}\rightarrow \mathcal{A}$ is called an {\it additive $\ast$-derivation} if
$D(XY)=D(X)Y+XD(Y)$ and  $D(X^\ast)=D(X)^\ast$ hold for all $X,Y\in\mathcal{A}.$ An additive map $G:\mathcal{A}\rightarrow \mathcal{A}$ is said to be an {\it additive generalized $\ast$-derivation} if there exists an additive $\ast$-derivation $D$ of $\mathcal{A}$ such that $G(XY)=G(X)Y+XD(Y)$ and  $G(X^\ast)=G(X)^\ast$ hold for all
$X,Y\in\mathcal{A}.$

For any \(X_1,X_2,\ldots,X_n\in\mathcal{A}\), define the following sequence of polynomials \(p_n\) recursively by
\[
p_1(X_1)=X_1
\]
and
\[
p_n(X_1,X_2,\ldots,X_n)
=[p_{n-1}(X_1,X_2,\ldots,X_{n-1}),X_n]_{\ast},
\]
for every integer \(n\geq 2\).

A map \(D:\mathcal{A}\to\mathcal{A}\), not necessarily additive or linear, is called a {\it nonlinear skew Lie \(n\)-derivation} (\(n\geq 2\)) if
\begin{equation}\label{eqn1.1}
D(p_n(X_1,X_2,\ldots,X_n))
=\sum_{i=1}^{n}
p_n(X_1,\ldots,X_{i-1},D(X_i),X_{i+1},\ldots,X_n)
\end{equation}
holds for all \(X_1,X_2,\ldots,X_n\in\mathcal{A}\).
In particular, a nonlinear skew Lie \(2\)-derivation is called a {\it nonlinear skew Lie derivation}, whereas a nonlinear skew Lie \(3\)-derivation is referred to as a {\it nonlinear skew Lie triple derivation}. These maps, together with nonlinear skew Lie \(n\)-derivations for \(n\geq 4\), are collectively referred to as {\it nonlinear skew Lie-type derivations}. 

A map \(G:\mathcal{A}\to\mathcal{A}\), not necessarily additive or linear, is called a {\it nonlinear generalized skew Lie \(n\)-derivation} if there exists a nonlinear skew Lie \(n\)-derivation \(D:\mathcal{A}\to\mathcal{A}\) such that
\begin{equation}\label{eqn1.2}
\begin{aligned}
G(p_n(X_1,X_2,\ldots,X_n))
={}&p_n(G(X_1),X_2,\ldots,X_n)+\sum_{i=2}^{n}
p_n(X_1,\ldots,X_{i-1},D(X_i),X_{i+1},\ldots,X_n)
\end{aligned}
\end{equation}
for all \(X_1,X_2,\ldots,X_n\in\mathcal{A}\). In particular, a nonlinear generalized skew Lie \(2\)-derivation is called a {\it nonlinear generalized skew Lie derivation}, whereas a nonlinear generalized skew Lie \(3\)-derivation is referred to as a {\it nonlinear generalized skew Lie triple derivation}. These maps, together with nonlinear generalized skew Lie \(n\)-derivations for \(n\geq 4\), are collectively referred to as {\it nonlinear generalized skew Lie-type derivations}. 
Clearly, every nonlinear skew Lie \(n\)-derivation is a particular example of a nonlinear generalized skew Lie \(n\)-derivation, by taking \(G=D\) in \eqref{eqn1.2}.

Over the past decade, the characterization of nonlinear skew Lie-type derivations on various classes of rings and algebras has attracted considerable attention (see, for instance, \cite{Ashraf2024,Madni2023,Jing2016,Lin1,Lin2,YZ12,LZC16} and the references therein). Yu and Zhang \cite{YZ12} established that every nonlinear skew Lie derivation on a factor von Neumann algebra is necessarily an additive \(\ast\)-derivation. This result was subsequently extended to nonlinear skew Lie triple derivations by Li et al. \cite{LZC16}.
Let \(\mathcal{H}\) be a complex Hilbert space, and let \(\mathcal{B}(\mathcal{H})\) denote the algebra of all bounded linear operators on \(\mathcal{H}\). In \cite{LLF14}, Li et al. proved that, if \(\mathcal{A}\subseteq\mathcal{B}(\mathcal{H})\) is a von Neumann algebra without central summands of type \(I_1\), then a map \(\delta:\mathcal{A}\to\mathcal{B}(\mathcal{H})\) is a nonlinear skew Lie derivation if and only if it is an additive \(\ast\)-derivation. Jing \cite{Jing2016} investigated nonlinear skew Lie derivations on standard operator algebras. More precisely, he showed that every nonlinear skew Lie derivation on a standard operator algebra \(\mathcal{A}\) on \(\mathcal{H}\), which is closed under the adjoint operation, is necessarily a linear \(\ast\)-derivation. This characterization was subsequently extended to nonlinear skew Lie-type derivations by Lin in \cite{Lin1}.
More recently, Madni et al. \cite{Madni2026} investigated nonlinear generalized skew Lie derivations on unital \(\ast\)-algebras with nontrivial projections. In particular, they established that every nonlinear generalized skew Lie derivation \(G\), with associated nonlinear skew Lie derivation \(D\), on a unital \(\ast\)-algebra \(\mathcal{A}\) admits the standard form
$G(X)=ZX+D(X),$ where \(Z\in\mathcal{Z}(\mathcal{A})\).

Motivated by the aforementioned developments, the present paper is devoted to characterizing nonlinear generalized skew Lie-type derivations on unital \(\ast\)-algebras. More precisely, we prove that, under certain mild assumptions, every nonlinear generalized skew Lie \(n\)-derivation \(G\) on a unital \(\ast\)-algebra \(\mathcal{A}\) containing a nontrivial projection admits the standard form
\[
G(X)=ZX+D(X),\qquad X\in\mathcal{A},
\]
where \(Z\in\mathcal{Z}(\mathcal{A})\) and \(D:\mathcal{A}\to\mathcal{A}\) is an additive \(\ast\)-derivation (Theorem \ref{thm3.1}). As applications, we obtain characterizations of nonlinear generalized skew Lie \(n\)-derivations on several important classes of unital \(\ast\)-algebras, including unital prime \(\ast\)-algebras (Corollary \ref{cor3.1}), standard operator algebras  (Corollary \ref{cor3.2}), factor von Neumann algebras (Corollary \ref{cor3.3}), and von Neumann algebras without central summands of type \(I_1\) (Corollary \ref{cor3.4}).


To establish our main result, we first reduce the problem of characterizing nonlinear generalized skew Lie \(n\)-derivations to that of characterizing a particular class of nonlinear maps. Indeed, let \(G:\mathcal{A}\to\mathcal{A}\) be a nonlinear generalized skew Lie \(n\)-derivation with associated nonlinear skew Lie \(n\)-derivation \(D:\mathcal{A}\to\mathcal{A}\). Then, in view of \eqref{eqn1.1} and \eqref{eqn1.2}, the map \(\Phi=G-D\) satisfies
\begin{equation}\label{eqn1.3}
\Phi(p_n(X_1,X_2,\ldots,X_n))
=p_n(\Phi(X_1),X_2,\ldots,X_n)
\end{equation}
for all \(X_1,X_2,\ldots,X_n\in\mathcal{A}\). Consequently, the characterization of nonlinear generalized skew Lie \(n\)-derivations reduces to the study of maps satisfying \eqref{eqn1.3}.
Accordingly, we call a map \(\Phi:\mathcal{A}\to\mathcal{A}\) satisfying \eqref{eqn1.3} a {\it nonlinear skew Lie \(n\)-centralizer}. In particular, a nonlinear skew Lie \(2\)-centralizer is called a {\it nonlinear skew Lie centralizer}, whereas a nonlinear skew Lie \(3\)-centralizer is referred to as a {\it nonlinear skew Lie triple centralizer}. 
These maps, together with nonlinear skew Lie \(n\)-centralizers for \(n\geq 4\), are collectively referred to as {\it nonlinear skew Lie-type centralizers}. 
Recently, nonlinear skew Lie centralizers and nonlinear skew Lie triple centralizers have been investigated in \cite{Madni2026} and \cite{LiLiWang2024}, respectively.

Motivated by these developments, in Section 2, we show that every nonlinear skew Lie \(n\)-centralizer \(\Phi\) on a unital \(\ast\)-algebra \(\mathcal{A}\) is of the form
$\Phi(X)=ZX,~ X\in\mathcal{A},$ where \(Z=Z^{\ast}\in\mathcal{Z}(\mathcal{A})\). In particular, such a map is a {\it linear \(\ast\)-centralizer}, that is, a linear map \(\Phi:\mathcal{A}\to\mathcal{A}\) satisfying $\Phi(XY)=\Phi(X)Y=X\Phi(Y)$ and $\Phi(X^{\ast})=\Phi(X)^{\ast}$ for all \(X,Y\in\mathcal{A}\).
Section 3 is devoted to the proof of our main result, Theorem \ref{thm3.1}, which provides the desired characterization of nonlinear generalized skew Lie \(n\)-derivations. Finally, in Section 4, we apply the main result to several well-known classes of unital \(\ast\)-algebras.

We close this section with the following elementary lemma which will be frequently used in subsequent discussions.
\begin{lemma}\label{lem1.1}
  For any $X\in\mathcal{A}$ and any integer $n\geq 2,$ we have
  \begin{itemize}
    \item[$(i)$] $p_n(X, I, I,\ldots, I)=2^{n-2}(X-X^{\ast});$ 
    \item[$(ii)$] $p_n(X, \frac{I}{2}, \frac{I}{2},\ldots, \frac{I}{2})=\frac{1}{2}(X-X^{\ast}).$
  \end{itemize}
\end{lemma}

\section{Nonlinear skew Lie-type centralizers}
In this section, we obtain a description of nonlinear skew Lie-type centralizers unital $\ast$-algebras. 
The main result of this section states as follows:
\begin{theorem}\label{prop2.1}
Let $\mathcal{A}$ be an unital $\ast$-algebra over the complex field $\mathbb{C}$. If a map $\Phi:\mathcal{A}\rightarrow\mathcal{A}$ satisfies
$\Phi(p_n(X_1,X_2, \ldots, X_n))=p_n(\Phi(X_1),X_2, \ldots, X_n)$ for all $X_1, X_2, \ldots, X_n\in\mathcal{A},$ then $\Phi$ is of the form  $\Phi(A)=Z A$ for all $A\in\mathcal{A},$ where $Z=Z^{\ast}\in\mathcal{Z}(\mathcal{A}).$ In particular, $\Phi$ is a linear $\ast$-centralizer. 
\end{theorem}
We will complete the proof by several lemmas.
\begin{lemma}\label{lem2.1}
  $\Phi(0)=0$
\end{lemma}
\begin{proof}
  $\Phi(0)=\Phi(p_n(0,0, \ldots, 0))=p_n(\Phi(0),0, \ldots, 0)=0.$
\end{proof}

\begin{lemma}\label{lem2.2}
For every \(\lambda\in\mathbb{R}\), we have $\Phi(\lambda I)\in\mathcal{Z}(\mathcal{A}).$
\end{lemma}

\begin{proof}
Let \(\lambda\in\mathbb{R}\) be arbitrary. Since $p_n(\lambda I,I,I,\ldots,I)=0,$
it follows from Lemmas \ref{lem1.1} and \ref{lem2.1} that
\begin{align*}
0
&=\Phi(p_n(\lambda I,I,I,\ldots,I))\\
&=p_n(\Phi(\lambda I),I,I,\ldots,I)\\
&=2^{n-2}\bigl(\Phi(\lambda I)-\Phi(\lambda I)^{\ast}\bigr).
\end{align*}
Consequently,
\begin{equation}\label{eq:2.1}
\Phi(\lambda I)=\Phi(\lambda I)^{\ast}.
\end{equation}

We now distinguish two cases.\\

{\bf Case 1.} Suppose that \(n=2\). For any \(A\in\mathcal{A}\), we have $p_2(\lambda I,A)=[\lambda I,A]_{\ast}=0.$ Therefore, using \eqref{eq:2.1}, we obtain
\[
0=\Phi([\lambda I,A]_{\ast})=[\Phi(\lambda I),A]_{\ast}=\Phi(\lambda I)A-A\Phi(\lambda I)^{\ast}=\Phi(\lambda I)A-A\Phi(\lambda I).
\]
Thus, \(\Phi(\lambda I)\) commutes with every element of \(\mathcal{A}\), and hence
\[
\Phi(\lambda I)\in\mathcal{Z}(\mathcal{A}).
\]

{\bf Case 2.} Suppose that \(n\geq 3\). Since $p_n(\lambda I,A,I,\ldots,I)=0$ for every \(A\in\mathcal{A}\), invoking \eqref{eq:2.1}, Lemmas \ref{lem1.1} and \ref{lem2.1}, we obtain
\begin{align*}
0
&=\Phi(p_n(\lambda I,A,I,\ldots,I))\\
&=p_n(\Phi(\lambda I),A,I,\ldots,I)\\
&=p_{n-1}([\Phi(\lambda I),A]_{\ast},I,\ldots,I)\\
&=2^{n-3}\bigl([\Phi(\lambda I),A]_{\ast}-([\Phi(\lambda I),A]_{\ast})^{\ast}\bigr)\\
&=2^{n-3}\bigl(\Phi(\lambda I)A-A\Phi(\lambda I)^{\ast}-A^{\ast}\Phi(\lambda I)^{\ast}
+\Phi(\lambda I)A^{\ast}\bigr)\\
&=2^{n-3}\bigl(\Phi(\lambda I)(A+A^{\ast})
-(A+A^{\ast})\Phi(\lambda I)\bigr).
\end{align*}
Hence, we conclude that
\[
\Phi(\lambda I)(A+A^{\ast})
=(A+A^{\ast})\Phi(\lambda I) \ \ \text{for every} \ \ A\in\mathcal{A}.
\]
In particular, taking \(A=A^{\ast}\), we see that
$
\Phi(\lambda I)A=A\Phi(\lambda I)
$
for every self-adjoint element \(A\in\mathcal{A}\).
Now, for an arbitrary \(A\in\mathcal{A}\), we can write
$
A=A_1+iA_2,
$
where
$
A_1=\frac{A+A^{\ast}}{2}$ and $A_2=\frac{A-A^{\ast}}{2i}
$
are both self-adjoint elements of \(\mathcal{A}\). Hence,
\begin{align*}
\Phi(\lambda I)A
&=\Phi(\lambda I)(A_1+iA_2)\\
&=A_1\Phi(\lambda I)+iA_2\Phi(\lambda I)\\
&=(A_1+iA_2)\Phi(\lambda I)\\
&=A\Phi(\lambda I).
\end{align*}
Thus, \(\Phi(\lambda I)\) commutes with every element of \(\mathcal{A}\), which yields
$
\Phi(\lambda I)\in\mathcal{Z}(\mathcal{A}).
$

Combining the two cases, we conclude that
\[
\Phi(\lambda I)\in\mathcal{Z}(\mathcal{A})\ \ \text{for every}\ \ \lambda\in\mathbb{R}.
\]
\end{proof}

\begin{lemma}\label{lem2.3}
We have
$
\Phi\left(-\frac{1}{2}iI\right)^{\ast}
=-\Phi\left(-\frac{1}{2}iI\right)
$
and
$
\Phi\left(-\frac{1}{2}iI\right)
=i\Phi\left(-\frac{1}{2}I\right).
$
Moreover,
$
\Phi\left(-\frac{1}{2}iI\right)\in\mathcal{Z}(\mathcal{A}).
$
\end{lemma}

\begin{proof}
Since
$
p_n\left(-\frac{1}{2}iI,\frac{1}{2}I,
\frac{1}{2}I,\ldots,\frac{1}{2}I\right)
=-\frac{1}{2}iI,
$
it follows from Lemma \ref{lem1.1} that
\begin{align*}
\Phi\left(-\frac{1}{2}iI\right)
&=\Phi\left(p_n\left(-\frac{1}{2}iI,
\frac{1}{2}I,\ldots,\frac{1}{2}I\right)\right)\\
&=p_n\left(\Phi\left(-\frac{1}{2}iI\right),
\frac{1}{2}I,\ldots,\frac{1}{2}I\right)\\
&=\frac{1}{2}\left(
\Phi\left(-\frac{1}{2}iI\right)
-\Phi\left(-\frac{1}{2}iI\right)^{\ast}
\right).
\end{align*}
Consequently,
\begin{equation}\label{eq:2.2}
\Phi\left(-\frac{1}{2}iI\right)^{\ast}
=-\Phi\left(-\frac{1}{2}iI\right).
\end{equation}

We now establish the second assertion by considering two cases.

{\bf Case 1.} Suppose that \(n=2\). Observe that
$p_2\left(-\frac{1}{2}iI,-\frac{1}{2}iI\right)=-\frac{1}{2}I.$
Thus, by \eqref{eq:2.2}, we obtain
\begin{align*}
\Phi\left(-\frac{1}{2}I\right)
&=\Phi\left(p_2\left(-\frac{1}{2}iI,
-\frac{1}{2}iI\right)\right)\\
&=p_2\left(\Phi\left(-\frac{1}{2}iI\right),
-\frac{1}{2}iI\right)\\
&=\left[\Phi\left(-\frac{1}{2}iI\right),
-\frac{1}{2}iI\right]_{\ast}\\
&=-\frac{1}{2}i\left(
\Phi\left(-\frac{1}{2}iI\right)
-\Phi\left(-\frac{1}{2}iI\right)^{\ast}
\right)\\
&=-i\Phi\left(-\frac{1}{2}iI\right).
\end{align*}

{\bf Case 2.} Suppose that \(n\geq 3\). By the definition of \(p_n\) and Lemma \ref{lem1.1}, we have
\begin{align*}
p_n\left(-\frac{1}{2}iI,\frac{1}{2}I,
\ldots,\frac{1}{2}I,-\frac{1}{2}iI\right)
&=
\left[p_{n-1}\left(-\frac{1}{2}iI,
\frac{1}{2}I,\ldots,\frac{1}{2}I\right),
-\frac{1}{2}iI\right]_{\ast}\\
&=
\left[-\frac{1}{2}iI,-\frac{1}{2}iI\right]_{\ast}
=-\frac{1}{2}I.
\end{align*}
Hence, again by Lemma \ref{lem1.1} and \eqref{eq:2.2}, we obtain
\begin{align*}
\Phi\left(-\frac{1}{2}I\right)
&=\Phi\left(p_n\left(-\frac{1}{2}iI,
\frac{1}{2}I,\ldots,\frac{1}{2}I,
-\frac{1}{2}iI\right)\right)\\
&=p_n\left(\Phi\left(-\frac{1}{2}iI\right),
\frac{1}{2}I,\ldots,\frac{1}{2}I,
-\frac{1}{2}iI\right)\\
&=\left[p_{n-1}\left(\Phi\left(-\frac{1}{2}iI\right),
\frac{1}{2}I,\ldots,\frac{1}{2}I\right),
-\frac{1}{2}iI\right]_{\ast}\\
&=\left[\Phi\left(-\frac{1}{2}iI\right),
-\frac{1}{2}iI\right]_{\ast}\\
&=-i\Phi\left(-\frac{1}{2}iI\right).
\end{align*}

Combining the two cases, we conclude that
\begin{equation}\label{eq:2.3}
\Phi\left(-\frac{1}{2}iI\right)
=i\Phi\left(-\frac{1}{2}I\right).
\end{equation}
Finally, by Lemma \ref{lem2.2}, we know that
$
\Phi\left(-\frac{1}{2}I\right)
\in\mathcal{Z}(\mathcal{A}).
$
Since the center \(\mathcal{Z}(\mathcal{A})\) is closed under scalar multiplication, it follows that
\begin{equation}\label{eq:2.4}
\Phi\left(-\frac{1}{2}iI\right)
=i\Phi\left(-\frac{1}{2}I\right)
\in\mathcal{Z}(\mathcal{A}).
\end{equation}

\end{proof}

\begin{lemma}\label{lem2.4}
For every \(A\in\mathcal{A}\), we have
\[
\Phi(iA)=i\Phi(I)A=iA\Phi(I)
\]
and
\[
\Phi(A)=\Phi(I)A=A\Phi(I).
\]
\end{lemma}

\begin{proof}
Let \(A\in\mathcal{A}\) be arbitrary. We first show that
\begin{equation}\label{eq:2.5}
\Phi(-iA)=2i\Phi\left(-\frac{1}{2}I\right)A.
\end{equation}

Suppose first that \(n=2\). Since
$
p_2\left(-\frac{1}{2}iI,A\right)
=\left[-\frac{1}{2}iI,A\right]_{\ast}
=-iA,
$
it follows from Lemma \ref{lem2.3} that
\begin{align*}
\Phi(-iA)
&=\Phi\left(p_2\left(-\frac{1}{2}iI,A\right)\right)\\
&=p_2\left(\Phi\left(-\frac{1}{2}iI\right),A\right)\\
&=\left[\Phi\left(-\frac{1}{2}iI\right),A\right]_{\ast}\\
&=\Phi\left(-\frac{1}{2}iI\right)A
-A\Phi\left(-\frac{1}{2}iI\right)^{\ast}\\
&=i\Phi\left(-\frac{1}{2}I\right)A
+iA\Phi\left(-\frac{1}{2}I\right)\\
&=2i\Phi\left(-\frac{1}{2}I\right)A,
\end{align*}
where the last equality follows from Lemma \ref{lem2.2}.
Now, suppose that \(n\geq 3\). By the definition of \(p_n\) and Lemma \ref{lem1.1}, we have
\begin{align*}
&p_n\left(-\frac{1}{2}iI,\frac{1}{2}I,\ldots,
\frac{1}{2}I,A\right)\\
&\qquad=
\left[p_{n-1}\left(-\frac{1}{2}iI,
\frac{1}{2}I,\ldots,\frac{1}{2}I\right),A\right]_{\ast}\\
&\qquad=\left[-\frac{1}{2}iI,A\right]_{\ast}
=-iA.
\end{align*}
Consequently, using Lemma \ref{lem1.1} and \eqref{eq:2.2}, we obtain
\begin{align*}
\Phi(-iA)
&=\Phi\left(p_n\left(-\frac{1}{2}iI,
\frac{1}{2}I,\ldots,\frac{1}{2}I,A\right)\right)\\
&=p_n\left(\Phi\left(-\frac{1}{2}iI\right),
\frac{1}{2}I,\ldots,\frac{1}{2}I,A\right)\\
&=\left[p_{n-1}\left(\Phi\left(-\frac{1}{2}iI\right),
\frac{1}{2}I,\ldots,\frac{1}{2}I\right),A\right]_{\ast}\\
&=\left[\Phi\left(-\frac{1}{2}iI\right),A\right]_{\ast}\\
&=2i\Phi\left(-\frac{1}{2}I\right)A.
\end{align*}
Thus, \eqref{eq:2.5} holds for every \(n\geq 2\).
Taking \(A=iI\) in \eqref{eq:2.5}, we obtain
\[
\Phi(I)=-2\Phi\left(-\frac{1}{2}I\right).
\]
Substituting this identity into \eqref{eq:2.5}, we get
\[
\Phi(-iA)=-i\Phi(I)A=-iA\Phi(I),
\]
where the second equality follows from Lemma \ref{lem2.2}.
Finally, replacing \(A\) by \(-A\) and \(iA\), respectively, in the above identity, we arrive at
\[
\Phi(iA)=i\Phi(I)A=iA\Phi(I)
\]
and
\[
\Phi(A)=\Phi(I)A=A\Phi(I),
\]
respectively.
\end{proof}

\begin{proof}[Proof of Theorem \ref{prop2.1}]
Let $\mathcal{A}$ be a unital $\ast$-algebra over the complex field $\mathbb{C}.$
Let $\Phi:\mathcal{A}\rightarrow\mathcal{A}$ be a map satisfying $\Phi(p_n(X_1,X_2, \ldots, X_n))=p_n(\Phi(X_1),X_2, \ldots, X_n)$ for all $X_1, X_2, \ldots, X_n\in\mathcal{A}.$ 
Then, it follows from Lemma \ref{lem2.4} that $\Phi(A)=\Phi(I)A=A\Phi(I)$ for all $A\in\mathcal{A}.$ Set $Z=\Phi(I).$ 
Then, by Lemma \ref{lem2.2}, $Z=Z^{\ast}\in\mathcal{Z(A)}$ and $\Phi(A)=ZA$ for all $A\in\mathcal{A}.$
Moreover, it is straightforward to prove that $\Phi$ is a linear $\ast$-centralizer. 
\end{proof}

As a direct consequence of Theorem \ref{prop2.1}, we have:
\begin{corollary}\cite[Theorem 3.1]{Madni2026}\label{cor2.1}
Let $\mathcal{A}$ be an unital $\ast$-algebra over the complex field $\mathbb{C}$. If a map $\Phi:\mathcal{A}\rightarrow\mathcal{A}$ satisfies
$\Phi([X_1,X_2]_{\ast})=[\Phi(X_1),X_2]_{\ast}$ for all $X_1, X_2\in\mathcal{A},$ then $\Phi$ is of the form  $\Phi(A)=Z A$ for all $A\in\mathcal{A},$ where $Z=Z^{\ast}\in\mathcal{Z}(\mathcal{A}).$ In particular, $\Phi$ is a linear $\ast$-centralizer. 
\end{corollary}

\begin{corollary}\cite[Theorem 2.1]{LiLiWang2024}\label{cor2.2}
Let $\mathcal{A}$ be an unital $\ast$-algebra over the complex field $\mathbb{C}$. If a map $\Phi:\mathcal{A}\rightarrow\mathcal{A}$ satisfies
$\Phi([[X_1,X_2]_{\ast},X_3]_{\ast})=[[\Phi(X_1),X_2]_{\ast}],X_3]_{\ast}$ for all $X_1, X_2,X_3\in\mathcal{A},$ then $\Phi$ is of the form  $\Phi(A)=Z A$ for all $A\in\mathcal{A},$ where $Z=Z^{\ast}\in\mathcal{Z}(\mathcal{A}).$ In particular, $\Phi$ is a linear $\ast$-centralizer. 
\end{corollary}

\section{Nonlinear generalized skew Lie-type derivations}
The aim of this section is to establish a characterization of nonlinear generalized skew Lie \(n\)-derivations on unital \(\ast\)-algebras. 
The main result of the article states as follows:

\begin{theorem}\label{thm3.1}
Let \(\mathcal{A}\) be a unital \(\ast\)-algebra with identity element \(I\) and a nontrivial projection \(P\) satisfying the following conditions:
\[
X\mathcal{A}P=\{0\}\quad\Longrightarrow\quad X=0
\]
and
\[
X\mathcal{A}(I-P)=\{0\}\quad\Longrightarrow\quad X=0.
\]
Then every nonlinear generalized skew Lie \(n\)-derivation \(G:\mathcal{A}\to\mathcal{A}\) admits a representation of the form
$G(A)=ZA+D(A)$ for all $A\in\mathcal{A},$ where \(Z=Z^{\ast}\in\mathcal{Z}(\mathcal{A})\) and \(D:\mathcal{A}\to\mathcal{A}\) is an additive \(\ast\)-derivation. In particular, \(G\) is an additive generalized \(\ast\)-derivation.
\end{theorem}

To establish the theorem, we need the following couple of lemmas:
\begin{lemma}\label{prop2.2}
Let $\mathcal{A}$ be a unital $\ast$-algebra. Then every nonlinear generalized skew Lie $n$-derivation  $G:\mathcal{A}\rightarrow\mathcal{A}$ with associated nonlinear skew Lie $n$-derivation $D:\mathcal{A}\rightarrow\mathcal{A}$ is of the form $G(A)=ZA+D(A)$  for all $A\in\mathcal{A},$ where $Z=Z^{\ast}\in\mathcal{Z}(\mathcal{A}).$
\end{lemma}

\begin{proof}
Let $G:\mathcal{A}\to \mathcal{A}$ be a nonlinear generalized skew Lie $n$-derivation with associated nonlinear skew Lie $n$-derivation $D:\mathcal{A}\to \mathcal{A}.$ Then, in view of \eqref{eqn1.1} and \eqref{eqn1.2}, the map $\Phi=G-D$ satisfies
\begin{equation*}
  \Phi(p_n(X_1, X_2, \ldots, X_n))=p_n(\Phi(X_1), X_2, \ldots, X_n)
\end{equation*}
for all  $X_1, X_2, \ldots, X_n\in\mathcal{A}.$ Hence, it follows from Theorem \ref{prop2.1} that
$\Phi(A)=Z A$ for all $A\in\mathcal{A},$ where $Z=Z^{\ast}\in\mathcal{Z}(\mathcal{A}).$
Consequently, $G(A)=ZA+D(A)$ for all $A\in\mathcal{A},$ where $Z=Z^{\ast}\in\mathcal{Z}(\mathcal{A}).$
\end{proof}

\begin{lemma}\label{prop2.3}
 Let \(\mathcal{A}\) be a unital \(\ast\)-algebra with identity element \(I\) and a nontrivial projection \(P\) satisfying the following conditions:
\[
X\mathcal{A}P=\{0\}\quad\Longrightarrow\quad X=0
\]
and
\[
X\mathcal{A}(I-P)=\{0\}\quad\Longrightarrow\quad X=0.
\] 
 Then every nonlinear skew Lie $n$-derivation $D:\mathcal{A}\to\mathcal{A}$ is an additive $\ast$-derivation.
\end{lemma}
\begin{proof}
If \(n=2\), the desired conclusion follows directly from \cite[Theorem 2.1]{Madni2026}. Now, suppose that \(n\geq 3\). By \cite[Theorem 1]{Madni2023}, \(D\) is an additive \(\ast\)-derivation, provided that \(D\bigl(\frac{1}{2}iI\bigr)\) is self-adjoint.
We observe, however, that the self-adjointness condition on \(D\bigl(\frac{1}{2}iI\bigr)\) is redundant. Indeed, it can be established by adapting the argument used in the proof of Claim 2 in \cite{Lin2}. Consequently, \(D\) is an additive \(\ast\)-derivation, and the desired conclusion follows.
\end{proof}

\begin{proof}[Proof of Theorem \ref{thm3.1}]
Let $G:\mathcal{A}\to \mathcal{A}$ be a nonlinear generalized skew Lie $n$-derivation with associated nonlinear skew Lie $n$-derivation $D:\mathcal{A}\to \mathcal{A}.$
Then, it follow from Lemma \ref{prop2.2} that $G$ is of the form $G(A)=Z A+D(A)$ for all $A\in\mathcal{A},$
where $Z=Z^{\ast}\in \mathcal{Z}(\mathcal{A}).$ Since $D:\mathcal{A}\to \mathcal{A}$ is a nonlinear skew Lie $n$-derivation,
in view of Lemma \ref{prop2.3}, $D$ is an additive $\ast$-derivation.
It only remains to show that $G$ is an additive generalized $\ast$-derivation.
We assert that $G$ is an additive generalized $\ast$-derivation with associated additive $\ast$-derivation $D.$
Clearly, $G$ is additive. For any $X_1, X_2\in\mathcal{A},$ we have
\begin{eqnarray*}
G(X_1)X_2+X_1D(X_2)&=&(ZX_1+D(X_1))X_2+X_1D(X_2)\\
&=&(ZX_1)X_2+D(X_1)X_2+X_1D(X_2)\\
&=&Z(X_1X_2)+D(X_1X_2)\\
&=&G(X_1X_2).
\end{eqnarray*}
We now show that $G(X^\ast)=G(X)^\ast$ for all  $X\in\mathcal{A}.$ 
Define a map $\Phi=G-D.$ Then, $\Phi$ satisfies \eqref{eqn1.3} and hence, by  Theorem \ref{prop2.1}, 
$\Phi$ is a linear $\ast$-centralizer. Hence,  $\Phi(X^\ast)={\Phi(X)}^\ast$ for all $X\in\mathcal{A}.$ 
Thus, for any  $X\in\mathcal{A},$ we get
\[G(X^\ast)=\Phi(X^\ast)+D(X^\ast)=\Phi(X)^\ast+D(X)^\ast=(\Phi(X)+D(X))^\ast=G(X)^\ast.\]
Consequently, $G$ is an additive generalized $\ast$-derivation with associated additive $\ast$-derivation $D.$
\end{proof}

\section{Applications}
In this section, we apply Theorem \ref{thm3.1} to some classical examples of unital $\ast$-alegbras.
Let $\mathcal{A}$ be a unital $\ast$-algebra with identity $I$ and a nontrivial projection $P$ (that is, $P\neq 0, I$ and ${P}^2=P={P}^\ast$). Set $P_1=P$ and write $P_2=I-P_1$.
Assume that $\mathcal{A}$ satisfies the following property
\begin{equation}\label{eqn4.1}
  X\mathcal{A}P_i=\{0\} \Longrightarrow X=0, \ \  i=1,2.
\end{equation}
Some specific examples of unital $\ast$-algebras with nontrivial projections satisfying \eqref{eqn4.1} are unital prime $\ast$-algebras, factor von Neumann algebras, von Neumann algebras with no central summands of type $I_1$ (see \cite{DL14,LLF14} for details).

Recall that an algebra $\mathcal{A}$ is prime if for any $X_1, X_2\in\mathcal{A},$ $X_1\mathcal{A}X_2=\{0\}$ implies that  either $X_1=0$ or $X_2=0.$ 
It is easy to verify that every prime $\ast$-algebra with a nontrivial projection satisfies \eqref{eqn4.1}.
Therefore, as a direct consequence of Theorem \ref{thm3.1}, we have:

\begin{corollary}\label{cor3.1}
Let $\mathcal{A}$ be a unital prime $\ast$-algebra with a nontrivial projection.
Then every nonlinear generalized skew Lie $n$-derivation $G:\mathcal{A}\rightarrow\mathcal{A}$ is of the form
$G(A)=ZA+D(A)$  for all $A\in\mathcal{A},$ where $Z=Z^{\ast}\in\mathcal{Z}(\mathcal{A})$ and $D:\mathcal{A}\to\mathcal{A}$ is an additive $\ast$-derivation. 
In particular, $G:\mathcal{A}\to\mathcal{A}$ is an additive generalized $\ast$-derivation.
\end{corollary}

Let $\mathcal{H}$ be a complex Hilbert space and $\mathcal{B(H)}$ be the algebra of all bounded linear operators on $\mathcal{H}$.
Let $\mathcal{F(H)}\subseteq \mathcal{B(H)}$ denote the subalgebra of all bounded finite rank operators. 
A subalgebra $\mathcal{A}\subseteq \mathcal{B(H)}$ is called a standard operator algebra if it contains $\mathcal{F(H)}.$ 
If a standard operator algebra $\mathcal{A}$ is closed  under the adjoint operation and contains the identity operator $I$, 
then it is a unital prime $\ast$-algebra with a nontrivial projection and $\mathcal{Z(A)}=\mathbb{C}I$. 
Further, we know that if $\mathcal{A}$ is a standard operator algebra on an infinite dimensional Hilbert space $\mathcal{H}$ and $D:\mathcal{A}\rightarrow\mathcal{B(H)}$ is an additive $\ast$-derivation, then there exists $T\in\mathcal{B(H)}$ satisfying $T+T^{\ast}=0$ such that $D(A)=AT-TA$ for all $A\in\mathcal{A}$ (see \cite[Proof of Theorem 2.14]{Jing2016}).
Hence, in view of Corollary \ref{cor3.1} and the above fact, we have the following result:

\begin{corollary}\label{cor3.2}
Let $\mathcal{H}$ be an infinite dimensional complex Hilbert space
and $\mathcal{A}$ be a standard operator algebra on $\mathcal{H}$ containing the identity operator $I.$
Suppose that $\mathcal{A}$ is closed under the adjoint operation.
Then every nonlinear generalized skew Lie $n$-derivation $G:\mathcal{A}\to \mathcal{A}$ is of the form
$G(A)=\lambda A+(AT-TA)$  for all $A\in\mathcal{A},$ where $\lambda\in\mathbb{C}$ and $T\in\mathcal{B(H)}$ satisfying $T+T^{\ast}=0.$
\end{corollary}

A von Neumann algebra $\mathcal{A}$ is a weakly closed, self adjoint algebra of operators on a Hilbert space $\mathcal{H}$ containing the identity operator $I$. 
A von Neumann algebra $\mathcal{A}$ is a factor if its center is $\mathbb{C}I$.
It is well known that a factor von Neumann algebra is prime.
Hence, as an immediate consequence of Corollary \ref{cor3.1}, we get:

\begin{corollary}\label{cor3.3}
Let $\mathcal{A}$ be a factor von Neumann algebra with dim$(\mathcal{A})\geq 2.$
Then every nonlinear generalized skew Lie $n$-derivation $G:\mathcal{A}\rightarrow\mathcal{A}$ is of the form
$G(A)=ZA+D(A)$  for all $A\in\mathcal{A},$ where $Z=Z^{\ast}\in\mathcal{Z}(\mathcal{A})$ and $D:\mathcal{A}\to\mathcal{A}$ is an additive $\ast$-derivation. 
In particular, $G:\mathcal{A}\to\mathcal{A}$ is an additive generalized $\ast$-derivation.
\end{corollary}

Further, it is well-known that every von Neumann algebra with no central summands of type $I_1$ satisfies \eqref{eqn4.1}
(see \cite{DL14,LLF14} for details).
Therefore, applying Theorem \ref{thm3.1},  we obtain:

\begin{corollary}\label{cor3.4}
Let $\mathcal{A}$ be a  von Neumann algebra having no central summands of type $I_1.$
Then every nonlinear generalized skew Lie $n$-derivation $G:\mathcal{A}\rightarrow\mathcal{A}$ is of the form
$G(A)=ZA+D(A)$  for all $A\in\mathcal{A},$ where $Z=Z^{\ast}\in\mathcal{Z}(\mathcal{A})$ and $D:\mathcal{A}\to\mathcal{A}$ is an additive $\ast$-derivation. 
In particular, $G:\mathcal{A}\to\mathcal{A}$ is an additive generalized $\ast$-derivation.
\end{corollary}

\noindent
{\bf Conflict of Interests}: The authors have no conflicts of interest to declare that are relevant to the content of this article.\\
\noindent
{\bf Data Availability Statement}: There is no data set  associated with this manuscript.


\begin{thebibliography}{99}


\bibitem{Ashraf2022} M. Ashraf, M.S. Akhter, M.A. Ansari:  Nonlinear bi-skew Lie-type derivations on von Neumann algberas. 
Comm. Algebra {\bf {50}}(11) (2022),  4766--4780.


\bibitem{Ashraf2023} M. Ashraf, M.S. Akhter, M.A. Ansari: Nonlinear bi-skew Jordan-type derivations on factor von Neumann algebras.
Filomat {\bf {37}}(17) (2022),  5591--5599.


\bibitem{Ashraf2024} M. Ashraf, M.S. Akhter, M.A. Ansari:  Non-global nonlinear skew Lie $n$-derivations on $\ast$-algebras. 
Comm. Algebra {\bf {52}}(9) (2024),  3734--3749. https://doi.org/10.1080/00927872.2024.2328802


\bibitem{Bresar2000}  M. Bre\v{s}ar, M. Fo\v{s}ner: On rings with involution equipped with some new product. Publ. Math. Debrecen.
{\bf{57}}(1--2)  (2000), 121--134.

\bibitem{LiZhao2021}  C. Li, Y. Zhao and F. Zhao: Nonlinear $\ast$-Jordan-type derivations on $\ast$-algebras. 
Rocky Mountain Journal Of Mathematics. {\bf 51}(2) (2021), 601--612.

\bibitem{LiLiWang2024}  C. Li, X. Li and J. Wang: Nonlinear skew Lie triple centralizers (derivations) on $\ast$-algebras. 
Filomat {\bf 38}(18) (2024), 6413--6421. https://doi.org/10.2298/FIL2418413L

\bibitem{DL14} L. Dai, F. Lu:  Nonlinear maps preserving Jordan $\ast$-products. J. Math. Anal. Appl. {\bf 409} (2014), 180--188.



\bibitem{Jing2016}  W. Jing: Nonlinear $\ast$-Lie derivations of standard operator algebras. Quaestiones Mathematicae {\bf 39}(8) (2016), 1037-1046, DOI: 10.2989/16073606.2016.1247119





\bibitem{Lin1}   W.-H. Lin: Nonlinear $\ast$-Lie-type derivations on standard operator algebras. Acta Math. Hungar.  {\bf{154}}(2)  (2018), 480--500. https://doi.org/10.1007/s10474-017-0783-6

\bibitem{Lin2}   W.-H. Lin: Nonlinear $\ast$-Lie-type derivations on von Neumann algebras. Acta Math. Hungar.  {\bf{156}}(1)  (2018), 112--131. https://doi.org/10.1007/s10474-018-0803-1



\bibitem{LLF14}  C. Li, F. Lu, X. Fang: Nonlinear $\xi$-Jordan $\ast$-derivations on von Neumann algebras. Linear Multilinear Algebra {\bf 62} (2014), 466--473.


\bibitem{LZC16} C.-J. Li,  F.-F. Zhao,   Q.-Y. Chen:  Nonlinear skew Lie triple derivations between factors. Acta Math. Sin. (Engl. Ser.) {\bf{32}}(7) (2016), 821--830.

\bibitem{Madni2023} M.A. Madni,  A.S. Alali,   M.R. Mozumder:  Nonlinear skew Lie-type derivations on $\ast$-algebras. Mathematics 2023, 11, 3819. https://doi.org/10.3390/math11183819

\bibitem{Madni2026} M.A. Madni,  M.S. Akhtar,   M.R. Mozumder: Characterization of generalized skew Lie derivation on $\ast$-algebras. Bol. Soc. Paran. Mat. {\bf{44}} (2026), 1--13.


\bibitem{Molnar97} L. Moln{\'a}r:  Jordan $\ast$-derivation pairs on a complex $\ast$-algebra.
Aequationes Math. {\bf{54}}  (1997), 44--55.

\bibitem{MS97} L. Moln{\'a}r, P. \v{S}emrl:  Local Jordan $\ast$-derivations of standard operator algebras.
Proc. Amer. Math. Soc. {\bf{125}}  (1997), 447--454.

\bibitem{PS90} P. \v{S}emrl:  On Jordan $\ast$-derivations and an application. Colloq. Math. {\bf{59}} (1990), 241--251.

\bibitem{PS91} P. \v{S}emrl: Quadratic functionals and Jordan $\ast$-derivations. Studia Math. {\bf{97}}  (1991), 157--165.


\bibitem{PS93} P. \v{S}emrl:  Quadratic and quasi-quadratic functionals. Proc. Amer. Math. Soc. {\bf{119}} (1993), 1105--1113.

\bibitem{PS94} P. \v{S}emrl:  Jordan $\ast$-derivations of standard operator algebras. Proc. Amer. Math. Soc.{\bf{120}}(2) (1994), 515-519.





\bibitem{YZ12} W.-Y. Yu, J.-H. Zhang:  Nonlinear $\ast$-Lie derivations on factor von Neumann algebras.  Linear Algebra Appl. {\bf{437}}(8) (2012), 1979–1991,

\bibitem{Zhang16} F. Zhang: Nonlinear skew Jordan derivable maps on factor von neumann algebras. Linear Multilinear Algebra {\bf{64}}(10)  (2016), 2090-2103.





\end{thebibliography}
\end{document}